\documentclass[12pt]{article}

\usepackage{amsmath}    
\usepackage{amssymb}    
\usepackage{amsfonts}   
\usepackage{enumerate}
\usepackage{vmargin}

\def \qed {\hfill \vrule height6pt  width6pt depth0pt}
\renewcommand{\span}{\mathrm{span}}

\def \tens{\otimes}
\newcommand{\m}{\mathcal}
\def \bl{\langle}
\def \br{\rangle}
\def \dim{{\mathrm{dim}} \, }
\def \ran{{\mathrm{ran}} \, }
\def \ker{{\mathrm{ker}} \, }
\def \C{{\mathbb C}}
\def \N{{\mathbb N}}

\def \R{\mathbb R}
\def \ksi{\xi}
\renewcommand{\le}{\ensuremath{\leqslant}}
\renewcommand{\ge}{\ensuremath{\geqslant}}

\title{Factoriality of $q$-Gaussian von Neumann algebras}
\author{\'Eric Ricard}
\newtheorem{thm}{Theorem}

\newtheorem{cor}[thm]{Corollary}

\newcounter{remark}
\newenvironment{rk}{{\noindent
\bf Remark \arabic{remark} \ }\addtocounter{remark}{1}}
{
\smallskip
}
\newenvironment{pf}[1][]{\noindent {\it Proof #1} : }{\hbox{~}\qed
\smallskip
}
\newcommand{\g}{\Gamma_q(\mathcal H_\R)}
\newcommand{\gr}{\Gamma_{q,r}(\mathcal H_\R)}

\renewcommand{\H}{\mathcal H}
\newcommand{\F}{\mathcal F_q(\mathcal H_\R)}
\date{}
\begin{document}

\maketitle

\begin{abstract}
We prove that the von Neumann algebras generated by $n$ $q$-Gaussian 
elements, are factors for $n\ge 2$.
\end{abstract}
\section{Introduction}

In the early 70's, Frish and Bourret considered operators  satisfying the
$q$-canonical commutation relations, for $-1<q<1$ : 
$$l(e)l^*(f)- q l^*(f)l(e)= (e,f) Id.$$
Nevertheless their existence was proved only 20 years later by 
Bo\.zejko and Speicher in \cite{BS2}. Since, many people studied
the von Neumann algebra $\g$, generated by $q$-Gaussian random
variables $\{l(e)+l^*(e) ; e\in\H_\R\}$,
 and some of their generalizations. It is well know that  
$\g$ is type $II_1$. 
One of the interesting point is that
these algebras realize a kind of interpolating scale 
between $\Gamma_1(\H)$ which
is commutative and $\Gamma_{-1}(H)$ the hyperfinite $II_1$ factor. For
$q=0$, we recover the algebra generated by Voiculescu's semicircular elements, 
which is a central object in the free probability theory.

Among the known results, Bo\.zejko and Speicher showed that $\g$ is non 
injective under some condition on the dimension of $\H$, which was removed
by Nou \cite{Nou}. 
 Recently, Shlyakhtenko \cite{Sh} proved that they are solid 
for some values of $q$. The question of the factoriality of $\g$ was studied
by  Bo\.zejko, K\"ummerer and Speicher \cite{BKS}, they showed that 
if $\H$ is infinite dimensional then $\g$ is a factor. This condition 
was partially released by \'Sniady \cite{Sn}, 
who showed that this is still true if
the dimension of $\H$ is greater than a function of $q$.

\section{Preliminaries}
In this paper, $-1<q<1$ is a fixed real number, we will use standard 
notation and refer to the papers \cite{BS,BKS,Nou} for general 
background.

 Let $\H$ be the  complexification of some real Hilbert space $\H_\R$.
By $\H^{\tens_2n}$ ($n\ge 1$), we denote the hilbertian 
$n$-tensor product of $\H$
with itself, this space is equipped with a scalar product that we write 
$(.,.)$.    Let $P_n: \H^{\tens_2n}\to \H^{\tens_2n}$ be given by 
$$P_n(e_1\tens ...\tens e_n)=\sum_{\sigma\in S_n} q^{|\sigma|} e_{\sigma(1)}
\tens...\tens  e_{\sigma(n)}=\sum_{\sigma\in S_n} q^{|\sigma|} \phi(\sigma)
( e_1\tens ...\tens e_n),$$
where $S_n$ is the symmetric group on $n$ elements, $|\sigma|$ is the number
of inversion of $\sigma$, and  $\phi$ is the natural action of 
$S_n$ on $\H^{\tens_2n}$. It was shown in \cite{BS}, that this operator 
is bounded and strictly positive, therefore we denote by 
$\H^{\tens_n}$, the hilbert space  $\H^{\tens_2n}$ equip with the new scalar
product $\bl .,.\br$ given by
$$\forall x,y\in \H^{\tens_n} \qquad   \bl x,y\br = ( x , P_n(y)).$$
 From now on, if $x\in\H^{\tens_n}$, $\|x\|$ is the norm of $x$ with respect 
to this new scalar product. For instance, if $e\in \H$ and $\|e\|=1$, then
$$\|e^{\tens n}\|^2=[n]_q!,$$
where $[k]_q=\frac {1-q^k}{1-q}$ and $[n]_q!=[1]_q...[n]_q$.

\begin{rk}
We will use as a key point that the sequence $([n]_q!)$ behave like a 
geometric sequence. 
\end{rk}

Moreover, it is known that the following algebraic relation holds :
$$P_n= R_{n,k} ( P_{n-k}\tens P_k) \quad \textrm{with } 
R_{n,k}=\sum_{\sigma \in S_{n}/S_{n-k}\times S_{k}}q^{|\sigma|}
\phi(\sigma^{-1}),$$
and the sum runs over the representatives of the right cosets of
$S_{n-k}\times S_{k}$ in $S_{n}$ with minimal number of inversions.  
As a consequence, since $\|R_{n,k}\|_{B(\H^{\tens_2n})}\le C_q=\prod_{i\ge 1}
(1-|q|^i)^{-1}$, we get that the formal identity map 
$$ Id : \H^{\tens n-k}\tens_2 \H^{\tens k} \to  \H^{\tens n}$$ 
has norm bounded by $\sqrt{C_q}$.

\begin{rk}\label{norm}
 As an application, we get that, if $e_1,...e_n$ and $e$
are  norm 1 vectors in $\H$, then 
$$\| e_1\tens...\tens e_n \tens e^{\tens m}\|_{\H^{\tens{n+m}}}\le C_q^{n/2}
\sqrt{[m]_q!}.$$ 
\end{rk}

The $q$-deformed Fock space is the Hilbert space defined by
$$\F = \C \Omega \oplus \oplus_{n\ge 1} \H^{\tens_n},$$
where $\Omega$ is a unital vector, considered as the vacuum.
 Vectors in $\H$ will be called letters and an elementary 
tensor of letters in $\H^{\tens n}$ will be called a word of length $n$.

For $e\in \H_\R$, we consider left and right creation operators on $\F$, given 
by :
\begin{eqnarray*}
l(e)(e_1\tens...\tens e_n)&=& e\tens e_1\tens...\tens e_n \\
l_r(e)(e_1\tens...\tens e_n)&=& e_1\tens...\tens e_n\tens e
\end{eqnarray*}
They are bounded endomorphisms of $\F$, more precisely if $\|e\|=1$ then
$$\|l_r(e)\|=\| l(e)\| = \left\{\begin{array}{ll}
1 & \textrm{if } q\le 0 \\
\frac 1 {\sqrt{1-q}} & \textrm{if } q\ge 0
\end{array}\right.$$
 Their adjoints in $B(\F)$ 
are the annihilation operators : 
\begin{eqnarray*}
l^*(e)(e_1\tens...\tens e_n)&=& \sum_{1\le i\le n}
q^{i-1}(e,e_i) \tens e_1\tens...\tens \hat{e_i}\tens  ..\tens e_n \\
l_r^*(e)(e_1\tens...\tens e_n)&=&\sum_{1\le i\le n}
q^{n-i}(e,e_i) \tens e_1\tens...\tens \hat{e_i}\tens  ..\tens e_n 
\end{eqnarray*}
 where $\hat{e_i}$ denotes a removed letter, if $n=0$, we put 
$l^*(e)\Omega=l_r^*(e)\Omega=0$.

 The operators $l(e)$ satisfy
the $q$-commutation relations :
$$l(e)l^*(f)- q l^*(f)l(e)= (e,f) Id$$

For $e\in \H_\R$, let
$$ W(e)=l(e)+l^*(e) \qquad \textrm{and}\qquad  W_r(e)=l_r(e)+l_r^*(e).$$
So for $e\in \H_\R$, $W(e)$ is self-adjoint.

 $\g$ stands for the von Neumann algebra generated by $(W(e))_{e\in \H_r}$
$$\g= \{ \; W(e) \; ;\; e\in \H_\R \;\}''.$$
 
And, $\gr$ stands for the von Neumann algebra generated 
by $(W_r(e))_{e\in \H_\R}$
$$\gr= \{ \; W_r(e) \; ;\; e\in \H_\R \;\}''.$$

 We recall some classical results on those algebras, 
\begin{itemize}
\item The commutant of $\g$ is $\g'=\gr$.
\item The vacuum vector $\Omega$ is separating and cyclic for both 
$\g$ and $\gr$.
\item The vector state $\tau(x)=\bl x\Omega, \Omega\br$ is a trace for both 
$\g$ and $\gr$. 
\end{itemize}

According to the  second point, any $x\in \g$ is uniquely determined
by $\ksi=x.\Omega\in \F$, so we will call it $x=W(\ksi)$
(and similarly for $\gr$, $x=W_r(\ksi)$). This notation is consistent with 
the definition of $W(e)=l(e)+l^*(e)$. The subspace 
$\g. \Omega \subset \F$ of all such $\ksi$ contains all tensors of 
finite rank, so it contains all words. 
If $e_1\tens...\tens e_n$ is a word in $\F$, there is a nice 
description of $W( e_1\tens...\tens e_n)$ in terms of $l(e_i)$ called
the Wick formula :  
 $$W(e_{1}\tens ...e_{n})= \sum_{m=0}^n \sum_{\sigma \in S_n/S_{n-m}\times
S_m} q^{|\sigma|}
l(e_{{\sigma(1)}})... l(e_{{\sigma(n-m)}})l^*(e_{{\sigma(n-m+1)}})
...l^*(e_{{\sigma(n)}}),$$
where $\sigma$ is the representative of the right coset of
$S_{n-m}\times S_{m}$ in $S_{n}$ with minimal number of inversions.
There is a similar formula for $W_r$.

 Actually, the algebras $\g$ and $\gr$ are in standard form in 
$B(\F)$, but we won't use it. 
If we denote by $S$,  the anti-symmetry that inverses the order 
of words in $\H_\R$, then for any $\ksi\in  \g. \Omega$ :
$$W(\ksi)^*=W(S\ksi)\qquad \textrm{and} \qquad S.W(\ksi).S=W_r(S\ksi).$$
In particular $ \g. \Omega= \gr. \Omega$.

\begin{rk}
For $\ksi,\eta\in \g.\Omega$, we will frequently use  
$$W(\ksi)\eta=W(\ksi)W_r(\eta)\Omega=W_r(\eta)W(\ksi)\Omega=W_r(\eta)\ksi.$$
\end{rk}

\bigskip

Let $T:\H_\R\to \H_\R$, be a $\R$-linear contraction, then there is
a canonical $\C$-linear contraction, $\m F_q(T)$, on $\F$  extending $T$, 
called the first quantization ; formally 
$$ \m F_q(T) = Id_{\C \Omega} \oplus \oplus_{n\ge 1} \tilde T^{\tens n}$$
with $\tilde T$, the complexification of $T$ on $\H$.

The second quantization of $T$, is the unique 
unital completely positive map $\Gamma_q(T)$ on 
$\g$ satisfying, for $\ksi \in \g.\Omega$
$$  \Gamma_q(T)(W(\ksi))=W(\m F_q(T)\ksi).$$

 For instance, if $\m K_\R\subset H_\R$, the second quantization associated
to the orthogonal projection $P_{\m K_\R}$ on $\m K_\R$ 
is a conditional expectation 
$$ \Gamma_q(P_{\m K_\R}): \g\to \Gamma_q(\m K_\R)=
\{W(e) \,;\, e\in \m K_\R\,\}''.$$

%

\section{The main result}

Let $e\in \H_\R$ of norm one and denote by $E_e$ the closed subspace 
of $\F$ spanned by the elements
$\{e^{\tens_n} \,;\, n\ge 0\}$, 
that is $E_e=\m F_q(\R e)$.  It is easy to check that for any 
$x=W(\ksi)\in W(e)''$, we have $\ksi \in E_e$. Conversely, assume 
$x=W(\ksi)$ and that 
$\ksi\in E_e$, then $x\in W(e)''$ : by the second quantization, we have 
a conditional expectation $\Gamma_q(P_{\R e}): \g\to W(e)''$, but then 
$$\Gamma_q(P_{\R e})(x).\Omega=\m F_q(P_{\R e}).\ksi=P_{E_e}.\ksi=
\ksi=x.\Omega,$$
as $\Omega$ is separating, $x=\Gamma_q(P_{\R e})(x)\in W(e)''$.

%

\begin{thm}
Assume that $\dim \H\ge 2$ and
let $e\in \H_\R$, $\|e\|=1$, then $W(e)''$ is a maximal abelian subalgebra
in $\g$.  
\end{thm}

\begin{cor}
$\g$ is a factor as soon as $\dim \H\ge 2$.
\end{cor}
\begin{pf}
Let $x\in \g\cap \g'$, then there is $\ksi\in \F$ such that
$x=W(\ksi)$. By the theorem, we must have $x\in W(e)''$ for every 
$e\in \H_\R$, 
but then $\ksi\in E_e$, so necessarily $x\in \C\Omega$. 
\end{pf}

\begin{pf}
Fix $(e_i)_{ i\ge 0}$ an orthonormal basis in $\H_\R$, with $e_0=e$.

Let $x=W(\ksi)\in \g\cap W(e)'$, we have to show that $\ksi\in E_e$.
For any $y=W(\eta)$ with $\eta \in E_e$, we have   
\begin{eqnarray*}
xy -yx =0 \\
(W(\ksi)W(\eta)-W(\eta)W(\ksi)).\Omega=0 \\
(W_r(\eta)-W(\eta)) \ksi=0
\end{eqnarray*}
So $\ksi\in \cap_{y=W(\eta)\in W(e)''} \ker 
(W_r(\eta)-W(\eta))$. By duality, we have to prove that
$$\overline{\span} \{ \ran  (W_r(\eta)-W(\eta)) \,;\, y=W(\eta)\in W(e)''\}
\supset E_e^\bot.$$

$E_e^\bot$ is the closed linear span of the set of elementary tensors 
$$F=\{ e_{i_1}\tens...\tens e_{i_n} ; n\ge 1, \textrm{ and } 
(i_1,...,i_n)\in \N^n\backslash \{(0,...,0)\}\}$$

Let $z=e_{i_1}\tens...\tens e_{i_n}$ 
be a word in $F$, it suffices to prove that $z$ is a weak-limit 
of elements in  
$\span \{ \ran  (W_r(\eta)-W(\eta)) \,;\, y=W(\eta)\in W(e)''\}$.

The von Neumann algebra $W(e)''$ is commutative and diffuse and 
separably generated (see \cite{BKS}), so we can 
assume that $W(e)''= L_\infty([0,1],dm)$, where $dm$ is the Lebesgue measure.
With this identification, the Rademacher functions $r_i$ belong to 
$W(e)''$, so we have $r_i=W(\eta_i)$ for some $\eta_i\in E_e$. Obviously
$W(\eta_i)$ is a self-adjoint symmetry and $W(\eta_i)^2=1$. 
Moreover, the sequence 
$(\eta_i)_{i\ge 1}$ converges to 0 for the weak topology on $\F$, since
$r_i$ is an orthormal basis in $L_2([0,1],dm)$.

Consider $$z_i=(W(\eta_i)-W_r(\eta_i))(W(\eta_i)(z)),$$
obviously $z_i \in 
{\span} \{ \ran  (W_r(\eta)-W(\eta)) \,;\, y=W(\eta)\in W(e)''\}$ and a simple 
calculation gives
$$z_i = W(\eta_i)^2 (z)- W_r(\eta_i)W(\eta_i)(z)=z- 
W_r(\eta_i)W(\eta_i)(z).$$
 We will show that $y_i= W_r(\eta_i).W(\eta_i)(z)$ tends
weakly to 0 in $\F$. As $\|y_i\|\le \|z\|$, it suffices to prove that for 
any word $t=e_{j_1}\tens...\tens e_{j_p}$, $\bl y_i, t \br\to 0$. We have,
\begin{eqnarray*}
\bl y_i,t\br&=& \bl W_r(\eta_i)).W(\eta_i)(z), t \br\\
            &=& \bl W_r(z)(\eta_i), W(t)(\eta_i)\br  
\end{eqnarray*}
  
This is the point where we use the Wick formula : 
$$W(e_{j_1}\tens ...\tens 
e_{j_n})= \sum_{m=0}^n \sum_{\sigma \in S_n/S_{n-m}\times
S_m} q^{|\sigma|}
l(e_{j_{\sigma(1)}})... l(e_{j_{\sigma(n-m)}})l^*(e_{j_{\sigma(n-m+1)}})
...l^*(e_{j_{\sigma(n)}})$$
and similarly for $z$. Since the number of terms appearing after developing
the sums is finite (it depends only on $n$ and $p$), we   
only need to show that
$$ I_i=\bl l_r(e_{i_1})...l_r(e_{i_m})
l_r^*(e_{i_{m+1}})...l_r^*(e_{i_n})(\eta_i),  l(e_{j_1})...l(e_{j_r})
l^*(e_{j_{r+1}})...l^*(e_{j_p})(\eta_i)\br \to 0, $$
as soon as at least one of the $i_k$'s is non zero. Let $v$ be the first 
$k$ such that $i_k\neq0$.
Since the letters in $\eta_i$ are only $e$, we can suppose that 
$v\le m$, otherwise $ l_r(e_{i_1})...l_r(e_{i_m})
l_r^*(e_{i_{m+1}})...l_r^*(e_{i_n})(\eta_i)=0$ (we have to cancel some 
$e_{i_v}$ in $\eta_i$ !). More generally, we can assume that
$e_{i_{m+1}}=...=e_{i_{n}}=e_{j_{r+1}}=...=e_{j_{p}}=e$.

Recall that 
$$l(e)^* e^{\tens n}=[n]_q  e^{\tens n-1}. $$
Now, we write $\eta_i=\sum_{k\ge 0} a_k^i e^{\tens k}$, interchanging the sums
and making simplifications
gives that ( with $a_{-n}=0$ if $n>0$. The $a_n$ are reals since $r_i$ is self
adjoint),
$$\begin{array}{cc} I_i=&\bl l_r(e_{i_1})...l_r(e_{i_m})
l_r^*(e_{i_{m+1}})...l_r^*(e_{i_n})(\eta_i), l(e_{j_1})...l(e_{j_r})
l^*(e_{j_{r+1}})...l^*(e_{j_p})(\eta_i)\br\\[12pt] = &
\displaystyle{\sum_{k\ge r,m}}
a_{k+n-2m}^i a_{k+p-2r}^i [k+n-2m]_q!/[k-m]_q![k+p-2r]_q!/[k-r]_q! \\ &.
\bl l_r(e_{i_1})...l_r(e_{i_m}) e^{\tens k-m}, 
l(e_{j_1})...l(e_{j_r}) e^{\tens k-r}\br\\[12pt]
=& \displaystyle{\sum_{k\ge r,m}}
a_{k+n-2m}^i a_{k+p-2r}^i [k+n-2m]_q!/[k-m]_q![k+p-2r]_q!/[k-r]_q!
\\& . \bl  l_r(e_{i_{v+1}})
...l_r(e_{i_m})e^{\tens k-m}, l_r^*(e_{i_v})...l_r^*(e_{i_1})
(e_{j_1}\tens...\tens e_{j_r}\tens e^{\tens k-r})\br
\end{array}$$

Assume that $k$ is big (say $k>N>2(n+p)$), by the  definition of $v$, we have
that $i_1=...=i_{v-1}=e$, so $$l_r^*(e_{i_{v-1}})...l_r^*(e_{i_1})
(e_{j_1}\tens...\tens e_{j_q}\tens e^{\tens k-r})$$ is obtained by cancelling 
$(v-1)$ times the letter $e$  in the word 
$e_{j_1}\tens...\tens e_{j_r}\tens e^{\tens k-r}$ using some 
 geometric weight $q^\alpha$,
$$\sum_{1\le {h_{v-1}} \le k-v-2}...
\sum_{1\le {h_2} \le k-1}\sum_{1\le {h_1} \le k}
\delta_{h_1,...} q^{(\sum h_i)-v+1}
(e_{j_1}\tens...\tens e_{j_r}\tens e^{\tens k-r})_{(h_1,...,h_{v-1})}$$
where $(e_{j_1}\tens...\tens e_{j_r}\tens e^{\tens k-r})_{(h_1,...,h_{v-1})}$
is obtained from $e_{j_1}\tens...\tens e_{j_r}\tens e^{\tens k-r}$ by 
removing the letter on the $h_1$-th position from the right, then the letter 
at $h_2$-th position in the remaining word and so on and where
$\delta_{h_1,...}$
is one if all the removed letters are $e$ and $0$ otherwise.
To have a non zero term in 
$$l_r^*(e_{i_v})
(e_{j_1}\tens...\tens e_{j_r}\tens e^{\tens k-r})_{(h_1,...,h_{v-1})}$$
we have to cancel a letter that is not an $e$, so it can happen only 
for the  terms coming from $e_{j_1}\tens...\tens e_{j_r}$ (if there 
are some left !), as this word of length 
$k-v+1$ ends with at least  $(k-r-v+1)$ $e$, we end up 
with a sum of at most $r$ words in front of which there is a factor less
than $|q|^{k-r-v+1}$. Moreover, by the remark \ref{norm}, the norm of such a 
word is less than $C_q^{r/2}\sqrt{[k-r-v+1]_q!}$.

If we sum up everything,  we get that 
$$\|l_r^*(e_{i_v})...l_r^*(e_{i_1})(e_{j_1}\tens...\tens e_{j_q}
\tens e^{\tens k-r})\| \le C(n,m,v,q) |q|^k \sqrt{[k]_q!}$$
where $C(n,m,v)$ does not depend on $k$ (because $[k]_q\le C_q$).

Now we can estimate $I_i$, by cutting the sum into two parts 
$A_i+B_i=\sum_{k\le N}|.|+ \sum_{k\ge N}|.|$.

Since $\eta_i\to 0$ weakly, each $a_{j}^i$ tends to 0,  then  
$$ A_i \le {\sum_{N>k\ge r,m}}
|a_{k+n-2m}^i| |a_{k+p-2r}^i| C(k,n,p)
\mathop{\to}^{i\to \infty} 0$$
and as $\|\eta_i\|\le 1$, we have $|a^i_k|\le 1/\sqrt{[k]_q!}$, so 
\begin{eqnarray*}
B_i&\le& \sum_{k\ge N} \frac {\sqrt{[k+n-2m]_q! [k+p-2r]_q!}}
{[k-m]_q!\,[k-r]_q!} 
\|l_r^*(e_{i_v})...l_r^*(e_{i_1})(e_{j_1}\tens...\tens e_{j_q}
\tens e^{\tens k-q})\|.\\ & &  \| l_r(e_{i_{v+1}})
...l_r(e_{i_m})e^{\tens k-m}\|\\ 
& \le &  \sum_{k\ge N} \frac {\sqrt{[k+n-2m]_q! [k+p-2r]_q!}}
{[k-m]_q!\,[k-r]_q!} C |q|^k \sqrt{[k]_q!} C(q)^m \sqrt{[k-m]_q!}\\
& \le&   \sum_{k\ge N}C |q|^k \le C |q|^N
\end{eqnarray*}

Consequently, we get that $\limsup |I_i| \le   C |q|^N$ for every $N$,
 so $I_i\to 0$.
\end{pf}
\nocite{*}
\bibliographystyle{plain}

\begin{flushleft}{ \'Eric Ricard}\\
 D\'epartement de Math\'ematiques de Besan{\c c}on
\\  Universit\'e de Franche-Comt\'e\\
 25030 Besan\c con cedex\\
 {\it e-mail}: {\tt eric.ricard@univ-fcomte.fr}\end{flushleft}
\end{document}